\documentclass[12pt]{amsart}
\usepackage{latexsym, amsmath, amscd, amssymb, amsthm} 
\usepackage[T2A]{fontenc}
 \usepackage[cp1251]{inputenc}
\usepackage[english]{babel}
\usepackage{color}
 \def\mathbi#1{\textbf{\em #1}}
\begin{document}

 \title{ The MAPLE package for $SL_2$-invariants and kernel of Weitzenb\"ock derivations.}

\author{Leonid Bedratyuk}  
\begin{abstract} 
  We  offer a Maple package {\tt SL\_2\_Inv\_Ker}  for calculating of minimal generating sets  for  the algebras of joint invariants/semi-invariants of binary forms   and  for calculations of  the kernels of Weitzenb\"ock  derivations.
\end{abstract}

\maketitle

\section{ Introduction}

Let $V_d \cong \langle v_0,v_1,\ldots,v_d \rangle $ be  $d+1$-dimensional  $SL_2$-module of binary forms of degree $d$ and let ${V_{\mathbi{d}} =
V_{d_1} \oplus V_{d_2} \oplus \cdots \oplus V_{d_n},}$ ${\mathbi{d}=(d_1,d_2,\ldots, d_n).}$ Denote by  $\mathcal{O}(V_{\mathbi{d}})^{SL_2}$   the algebra of polynomial $SL_2$-invariant functions on $V.$ Denote by  $\mathcal{O}(V_{\mathbi{d}})^{U_2}$   the algebra of polynomial $U_2$-invariant functions on $V,$ where $U_2$  is the maximal unipotent subgroup  of $SL_2.$  We  have  the  obvious inclusion $\mathcal{O}(V_{\mathbi{d}})^{SL_2} \subset \mathcal{O}(V_{\mathbi{d}})^{U_2}.$ Moreover $\mathcal{O}(V_{\mathbi{d}})^{U_2}\cong \mathcal{O}(V_1 \oplus V_{\mathbi{d}})^{SL_2}.$ The elements of the finitely generated algebras can be identified with the algebras of  joint  invariants  and  joint semi-invariants of the binary forms of degrees  $\mathbi{d}.$   

A linear derivation $\mathcal{D}$ of a  polynomial algebra  is called   Weitzenb\"ock derivation if its matrix (as linear map  on  the vector space  generated by  variables  of the  polynomial algebra)  is  nilpotent.
Any   Weitzenb\"ock derivation $\mathcal{D}$ is  completely determined by   the Jordan normal form of its  matrix. We will denote by  $\mathcal{D}_{\mathbi{d}},$  $\mathbi{d}:=(d_1,d_2,\ldots, d_n)$   the   Weitzenb\"ok derivation if  its  matrix consisting  of  $n$  Jordan blocks of size  $d_1+1,$ $d_2+1,$ $\ldots, d_n+1,$ respectively.

As is well known, the kernel of the derivations $\mathcal{D}_{\mathbi{d}},$  $\mathbi{d}:=(d_1,d_2,\ldots, d_n)$ is isomorphic to the algebra of joint semi-invariants (and  to the algebra of covariants) for    $n$ binary forms of degrees  $d_1,d_2,\ldots, d_n.$  Thus, the calculations of a minimal generated sets of the algebras $\mathcal{O}(V_{\mathbi{d}})^{SL_2}$  or $\mathcal{O}(V_{\mathbi{d}})^{U_2}$  can be reduced  to calculation of minimal generating set of the kernel of suitable Weitzenb\"ock derivation. For instance,  the algebra  $\mathcal{O}(V_d)^{U_2}$  can be identified with $\ker \mathcal{D},$ where the derivation $\mathcal{D}$   is defined by ${\mathcal{D}(v_i)=i v_{i-1}},$ $i=0,\ldots,d.$
 
 \section{ Algorithm}
 
  The algebra  $\mathcal{O}(V_{\mathbi{d}})^{U_2}$ is a finitely  generated multigraded algebra under the multidegree-order: 
$$
\mathcal{O}(V_{\mathbi{d}})^{U_2}=(\mathcal{O}(V_{\mathbi{d}})^{U_2})_{\mathbi{m},0}+(\mathcal{O}(V_{\mathbi{d}})^{U_2})_{\mathbi{m},1}+\cdots+(\mathcal{O}(V_{\mathbi{d}})^{U_2})_{\mathbi{m},j}+ \cdots,
$$
where  each subspace   $(\mathcal{O}(V_{\mathbi{d}})^{U_2})_{\mathbi{d},j}$ of joint semi-invariants of order $j$ and multidegree $\mathbi{m}:=(m_1,m_2,\ldots,m_n)$   is   finite-dimensional. The formal power series 
$$
\mathcal{P}(\mathcal{O}(V_{\mathbi{d}})^{U_2},z_1,z_2,\ldots,z_n,t)=\sum_{\mathbi{m},j=0}^{\infty }\dim((\mathcal{O}(V_{\mathbi{d}})^{U_2})_{\mathbi{m},j}) z_1^{m_1} z_2^{m_2}\cdots z_n^{m_n} t^j,
$$ 
is called the multivariariate Poincar\'e series   of the algebra of   join semi-variants  $\mathcal{C}_{\mathbi{d}}.$
Note that  each semi-invariant of order zero is an invariant.

Suppose that $\mathcal{P}(\mathcal{O}(V_{\mathbi{d}})^{SL_2},t)=\dfrac{P(t)}{Q(t)}$ and $P(t),$  $Q(t)$  are coprime polynomials. Denote  by $\beta_{\mathbi{d}}$ the degree  of the denominator $Q(t).$ 
Many experimental data leads  to the following conjecture:

{\bf Conjecture.} A degree upper bound for  irreducible invariants of the algebra $\mathcal{O}(V_{\mathbi{d}})^{SL_2}$ does not exceed $\beta_{\mathbi{d}}.$ 

\vspace{0.6 cm}

In the first version of the package we used the  following algorithm:

1. Calculate  multivariate Poincar\'e series of the algebras $\mathcal{O}(V_{\mathbi{d}})^{SL_2}$  ( or $\mathcal{O}(V_{\mathbi{d}})^{SL_2}$  or $\ker \mathcal{D}_{\mathbi{m}}.$)

2. For every  term  $z_1^{m_1} z_2^{m_2}\cdots z_n^{m_n} t^j$ of the Poincar\'e  series  calculate (by  linear algebra method) a basis of the vector  space of semi-invariants(or invariants  or elements of kernel ) and the multidegree $(m_1,m_2,\ldots,m_n).$ 

3. Separate irreducible polynomials.

4. Stop calculation  if  $m_1+m_2+\ldots + m_n >18.$

For  the  package procedures {\tt Min\_Gen\_Set\_Invariants\_S} we use also the following simplified algorithm:

1. Calculate   Poincar\'e series of the algebras $\mathcal{O}(V_{\mathbi{d}})^{SL_2}.$ 

2. For every  term  $z^m$ of the Poincar\'e series calculate  a basis of the vector  space of invariants of the degree $m.$

3. Separate irreducible polynomials.

4. Stop calculation  if  $m>18.$

The second algorithm works some fast  for  small values $d_i$ and $n>4.$

The package calculate the set  of irreducible invariants  up to degree $\min(18,\beta_{\mathbi{d}})$, but in all known computable  cases this set coincides with a minimal generating set, see,  for example, Brouwer's  webpage  
http://www.win.tue.nl/\!\text{$\sim$}\!\! aeb/math/invar/invarm.html

\section{ Installation.}  The package can be downloaded  from the web  http://sites.google.com/site/bedratyuklp/.
 
\begin{enumerate}
\item download the file {\tt SL\_2\_Inv\_Ker.mpl} and save it into your  Maple directory.
\item download the Xin's file (see link at the web page) {\tt Ell2.mpl} and save it into your  Maple directory.
 \item  run Maple
 
 \item  {\tt \textcolor{red}{ > restart:} \textcolor{red}{read "SL\_2\_Inv\_Ker.mpl":read "Ell2.mpl":}}

 \item  If necessary use {\tt \textcolor{red}{ > HelP();}}
\end{enumerate}

\section{Package Procedures  and Syntax}

\noindent
{\em Procedure name}: {\tt Min\_Gen\_Set\_Semi\_Invariants}\\
\noindent{\em Feature}: Computes  irreducible invariants  of the algebra $\mathcal{O}(V_{\mathbi{d}})^{U_2}$  up to degree $\min(18,\beta_{\mathbi{d}}).$\\
\noindent{\em Calling sequence}:
{\tt
Min\_Gen\_Set\_Semi\_Invariants$\left([d_1, d_2, \ldots,d_n] \right);$
}\\
\noindent{\em Parameters}:\\
{\raggedright
\begin{tabular}{lcl}
$ [d_1, d_2, \ldots,d_n]$ & - & a list of degrees of $n$  binary forms.\\
$ n $ & - & an integer, $n\leq 11.$\\
\end{tabular}\\
\vspace{5mm}
}

\noindent
{\em Procedure name}: {\tt Min\_Gen\_Set\_Invariants}\\
\noindent{\em Feature}: Computes  irreducible invariants  of the algebra $\mathcal{O}(V_{\mathbi{d}})^{SL_2}$  up to degree $\min(18,\beta_{\mathbi{d}}).$\\
\noindent{\em Calling sequence}:
{\tt
Min\_Gen\_Set\_Invariants$\left([d_1, d_2, \ldots,d_n] \right);$
}\\
\noindent{\em Parameters}:\\
{\raggedright
\begin{tabular}{lcl}
$ [d_1, d_2, \ldots,d_n]$ & - & a list of degrees of $n$  binary forms.\\
$ n $ & - & an integer, $n\leq 11.$\\
\end{tabular}\\
\vspace{5mm}
}

\noindent
{\em Procedure name}: {\tt Kernel\_LLN\_Der}\\
\noindent{\em Feature}: Computes  irreducible elements of the kernel of Weitzenb\"ock derivation   $\mathcal{D}_{\mathbi{d}}$  up to degree $\min(18,\beta_{\mathbi{d}}).$\\
\noindent{\em Calling sequence}:
{\tt
Kernel\_LLN\_Der$\left([d_1, d_2, \ldots,d_n] \right);$
}\\
\noindent{\em Parameters}:\\
{\raggedright
\begin{tabular}{lcl}
$ [d_1, d_2, \ldots,d_n]$ & - & a list of degrees of $n$  binary forms.\\
$ n $ & - & an integer, $n\leq 11.$\\
\end{tabular}\\
\vspace{5mm}
}

\noindent
{\em Procedure name}: {\tt Min\_Gen\_Set\_Invariants\_S}\\
\noindent{\em Feature}: By using the second algorithm the procedure computes  a set of   irreducible invariants  of the algebra $\mathcal{O}(V_{\mathbi{d}})^{SL_2}$  up to degree $\min(18,\beta_{\mathbi{d}}).$\\
\noindent{\em Calling sequence}:
{\tt
Min\_Gen\_Set\_Invariants\_S$\left([d_1, d_2, \ldots,d_n] \right);$
}\\
\noindent{\em Parameters}:\\
{\raggedright
\begin{tabular}{lcl}
$ [d_1, d_2, \ldots,d_n]$ & - & a list of degrees of $n$  binary forms.\\
$ n $ & - & an integer, $n\leq 11.$\\
\end{tabular}\\
\vspace{5mm}
}

\section{Examples}
\subsection{Compute $\mathcal{O}(V_4)^{SL_2}$} Use the command
\vspace{5mm}

{\tt \textcolor{red}{ > dd:=[4]:INV:=Min\_Gen\_Set\_Invariants(dd):}}\\

\begin{center}
\textcolor{blue}{         "calculating multivariate Poincare series...."}
\end{center}

                    \textcolor{blue} { "done!, upper  bound", 13}

  \textcolor{blue} { "---------------------------------------------------------------        ------------------"}

  \textcolor{blue} { "-----------------------------degree----------------------------
        --", 2}

         \textcolor{blue} {    " irreducible invariant of multidegree ", [2], "found"}

   \textcolor{blue} {"-----------------------------degree----------------------------
        --", 3}

          \textcolor{blue} {  " irreducible invariant of multidegree ", [3], "found"}

  \textcolor{blue} { "-----------------------------degree----------------------------
        --", 4}

  \textcolor{blue} { "-----------------------------degree----------------------------
        --", 5}

   \textcolor{blue} {"-----------------------------degree----------------------------
        --", 6}

   \textcolor{blue} {"-----------------------------degree----------------------------
        --", 7}

   \textcolor{blue} {"-----------------------------degree----------------------------
        --", 8}

   \textcolor{blue} {"-----------------------------degree----------------------------
        --", 9}

   \textcolor{blue} {"-----------------------------degree----------------------------
        --", 10}

   \textcolor{blue} {"-----------------------------degree----------------------------
        --", 11}

   \textcolor{blue} {"-----------------------------degree----------------------------
        --", 12}

   \textcolor{blue} {"-----------------------------degree----------------------------
        --", 13}

     \textcolor{blue} {"Total number of irreducible invariants   in     minimal generating set  ", 2}

\vspace{0.5cm}
To extract the two invariants  from the set INV one  should to use the commands

{\tt \textcolor{red} { >INV[2];}}

 \textcolor{blue} {$$\left\{ 6\,{x_{{2}}}^{2}+2\,x_{{0}}x_{{4}}-8\,x_{{1}}x_{{3}} \right\} $$}

{\tt  \textcolor{red} {>INV[3]};}  
 \textcolor{blue} {$$\left\{ -6\,{x_{{2}}}^{3}-6\,x_{{0}}{x_{{3}}}^{2}-6\,{x_{{1}}}^{2}x_{{4}}+6\,x_{{0}}x_{{2}}x_{{4}}+12\,x_{{1}}x_{{2}}x_{{3}} \right\} $$}

\subsection{Compute $\mathcal{O}(V_3 \oplus V_4)^{SL_2}$} Use the command
\vspace{5mm}

{\tt  \textcolor{red} {> dd:=[3,4]:INV:=Min\_Gen\_Set\_Invariants(dd):}}

{\begin{center}
\textcolor{blue} {"calculating multivariate Poincare series...."}

                   \textcolor{blue}{   "done!, upper  bound", 13}

 \textcolor{blue}{ "---------------------------------------------------------------------------------"}

  \textcolor{blue}{"-----------------------------degree------------------------------", 2}

       \textcolor{blue}{    "  irreducible invariant of multidegree  ", [0, 2], "found"}

  \textcolor{blue}{"-----------------------------degree------------------------------", 3}

       \textcolor{blue}{   " irreducible invariant of multidegree  ", [0, 3], "found"}

  \textcolor{blue}{"-----------------------------degree------------------------------", 4}

      \textcolor{blue}{     " irreducible invariant of multidegree", [4, 0], "found"}

  \textcolor{blue}{"-----------------------------degree------------------------------", 5}

          \textcolor{blue}{ " irreducible invariant of multidegree ", [4, 1], "found"}

          \textcolor{blue}{ " irreducible invariant of multidegree", [2, 3], "found"}

  \textcolor{blue}{"-----------------------------degree------------------------------", 6}

           \textcolor{blue}{ " irreducible invariant of multidegree", [4, 2], "found"}

            \textcolor{blue}{" irreducible invariant of multidegree", [4, 2], "found"}

 \textcolor{blue}{ "-----------------------------degree------------------------------", 7}

           \textcolor{blue}{ " irreducible invariant of multidegree", [4, 3], "found"}

           \textcolor{blue}{ " irreducible invariant of multidegree", [4, 3], "found"}

         \textcolor{blue}{   " irreducible invariant of multidegree", [4, 3], "found"}

 \textcolor{blue}{ "-----------------------------degree------------------------------", 8}

         \textcolor{blue}{ " irreducible invariant of multidegree ", [6, 2], "found"}

          \textcolor{blue}{  " irreducible invariant of multidegree", [4, 4], "found"}

           \textcolor{blue}{ "irreducible invariant of multidegree", [4, 4], "found"}

  \textcolor{blue}{"-----------------------------degree------------------------------", 9}

            \textcolor{blue}{" irreducible invariant of multidegree", [6, 3], "found"}

            \textcolor{blue}{" irreducible invariant of multidegree", [6, 3], "found"}

            \textcolor{blue}{" irreducible invariant of multidegree", [6, 3], "found"}

            \textcolor{blue}{" irreducible invariant of multidegree", [4, 5], "found"}

  \textcolor{blue}{"-----------------------------degree------------------------------", 10}

            \textcolor{blue}{" irreducible invariant of multidegree", [6, 4], "found"}

            \textcolor{blue}{" irreducible invariant of multidegree", [6, 4], "found"}

  \textcolor{blue}{"-----------------------------degree------------------------------", 11}

            \textcolor{blue}{" irreducible invariant of multidegree", [6, 5], "found"}

  \textcolor{blue}{"-----------------------------degree------------------------------", 12}

  \textcolor{blue}{"-----------------------------degree------------------------------", 13}
\end{center}

\begin{center}
\textcolor{blue} {"Total number of  invariants  in a minimal generating set ", 20}
\end{center}

To get the  number  of invariants of degree  $i$  use the command {\tt nops(INV[i])},  for  instance

{\tt  \textcolor{red}{> nops(INV[5])};}

\begin{center}
\textcolor{blue}{2}
\end{center}
 
 To extract and manipulate by these two  invariants use  the following commands
 
{\tt  \textcolor{red}{> Inv5\_1:=INV[5][1]:Inv5\_2:=INV[5][2];}}
\textcolor{blue}{
\begin{gather*}
Inv5\_2:=6\,{x_{{1}}}^{2}{x_{{3}}}^{2}y_{{0}}+18\,{x_{{1}}}^{2}{x_{{2}}}^{2}y_{
{2}}-12\,{x_{{1}}}^{3}x_{{3}}y_{{2}}-12\,{x_{{1}}}^{3}x_{{2}}y_{{3}}+6
\,{x_{{0}}}^{2}{x_{{2}}}^{2}y_{{4}}+6\,{x_{{2}}}^{4}y_{{0}}+\\+6\,{x_{{1}
}}^{4}y_{{4}}+12\,{x_{{1}}}^{2}x_{{2}}x_{{3}}y_{{1}}-12\,x_{{0}}x_{{1}
}{x_{{3}}}^{2}y_{{1}}+12\,x_{{0}}{x_{{2}}}^{2}x_{{3}}y_{{1}}+12\,x_{{0
}}{x_{{1}}}^{2}x_{{3}}y_{{3}}-12\,x_{{0}}{x_{{1}}}^{2}x_{{2}}y_{{4}}-\\-
12\,{x_{{0}}}^{2}x_{{2}}x_{{3}}y_{{3}}-12\,x_{{1}}{x_{{2}}}^{2}x_{{3}}
y_{{0}}+12\,x_{{0}}x_{{1}}{x_{{2}}}^{2}y_{{3}}-12\,x_{{0}}{x_{{2}}}^{3
}y_{{2}}+6\,{x_{{0}}}^{2}{x_{{3}}}^{2}y_{{2}}-12\,x_{{1}}{x_{{2}}}^{3}
y_{{1}}
\end{gather*}
}

 \subsection{Compute $\mathcal{O}(V_1 \oplus V_1 \oplus V_2)^{U_2}$} Use the command

 {\tt \textcolor{red}{> dd:=[1,1,2]:COV:=Min\_Gen\_Set\_Semi\_Invariants(dd):}}

\begin{center}
  \textcolor{blue}{"calculating multivariate Poincare series...."}

                     \textcolor{blue}{ "done!, upper  bound", 13}

  \textcolor{blue}{"---------------------------------------------------------------------------------"}

  \textcolor{blue}{"-----------------------------degree------------------------------", 2}

  \textcolor{blue}{"irreducible semi-invariant of multidegree", [0, 0, 2],
        "and order ", 0, "found"}

  \textcolor{blue}{"irreducible semi-invariant of multidegree", [1, 1, 0],
        "and order ", 0, "found"}

  \textcolor{blue}{"irreducible semi-invariant of multidegree", [0, 1, 1],
        "and order ", 1, "found"}

  \textcolor{blue}{"irreducible semi-invariant of multidegree", [1, 0, 1],
        "and order ", 1, "found"}

  \textcolor{blue}{"-----------------------------degree------------------------------", 3}

  \textcolor{blue}{"irreducible semi-invariant of multidegree", [2, 0, 1],
        "and order ", 0, "found"}

  \textcolor{blue}{"irreducible semi-invariant of multidegree", [1, 1, 1],
        "and order ", 0, "found"}

  \textcolor{blue}{"irreducible semi-invariant of multidegree", [0, 2, 1],
        "and order ", 0, "found"}

  \textcolor{blue}{"-----------------------------degree------------------------------", 4}

  \textcolor{blue}{"-----------------------------degree------------------------------", 5}

  \textcolor{blue}{"-----------------------------degree------------------------------", 6}

  \textcolor{blue}{"-----------------------------degree------------------------------", 7}

  \textcolor{blue}{"-----------------------------degree------------------------------", 8}

  \textcolor{blue}{"-----------------------------degree------------------------------", 9}

  \textcolor{blue}{"-----------------------------degree------------------------------", 10}

  \textcolor{blue}{"-----------------------------degree------------------------------", 11}

  \textcolor{blue}{"-----------------------------degree------------------------------", 12}

  \textcolor{blue}{"-----------------------------degree------------------------------", 13}

   \textcolor{blue}{  "number  of semi-invariant of  minimal generating set ", 10}
\end{center}
\vspace{0.5cm}
Below is this  minimal generating set listed by degree

\textcolor{red}{> COV[1];nops(\%);}

\begin{center}
\textcolor{blue}{$\{x_0,y_0,u_0 \}$}
\end{center}
\begin{center}
\textcolor{blue}{$3$}
\end{center}

\textcolor{red}{> COV[2];nops(\%);}

\begin{center}
\textcolor{blue}{$\left\{ -x_{{0}}y_{{1}}+x_{{1}}y_{{0}},-y_{{1}}u_{{0}}+y_{{0}}u_{{1}},-x_{{1}}u_{{0}}+x_{{0}}u_{{1}},2\,u_{{0}}u_{{2}}-2\,{u_{{1}}}^{2}
 \right\}
 $}
\end{center}
\begin{center}
\textcolor{blue}{$4$}
\end{center}

\textcolor{red}{> COV[3];nops(\%);}

\begin{center}
\textcolor{blue}{\begin{multline*}\{ -2\,y_{{0}}y_{{1}}u_{{1}}+{y_{{1}}}^{2}u_{{0}}+{y_{{0}}}^{2}u_{{2}},-x_{{1}}y_{{1}}u_{{0}}+x_{{1}}y_{{0}}u_{{1}}+x_{{0}}y_{{1}}u_{{
1}}-x_{{0}}y_{{0}}u_{{2}},\\{x_{{1}}}^{2}u_{{0}}+{x_{{0}}}^{2}u_{{2}}-2
\,x_{{0}}x_{{1}}u_{{1}} \}
\end{multline*}}
\end{center}
\begin{center}
\textcolor{blue}{$3$}
\end{center}

 \subsection{Compute the kernel  of the Weitzenb\"ock derivation $\mathcal{D}_{\mathbi{d}},$ $\mathbi{d}=(1,3)$}  Use the command
 
 {\tt \textcolor{red}{ > dd:=[1,3]:Ker:=Kernel\_LLN\_Der(dd):}}

\begin{center}

          \textcolor{blue}{  "calculating multivariate Poincare series...."}

                     \textcolor{blue}{ "done!, upper  bound", 13}

  \textcolor{blue}{"---------------------------------------------------------------
        ------------------"}

  \textcolor{blue}{"-----------------------------degree----------------------------", 2}

  \textcolor{blue}{"irreducible element of multidegree", [1, 1], "and order ",
        2, "found"}

  \textcolor{blue}{"irreducible element of multidegree", [0, 2], "and order ",
        2, "found"}

  \textcolor{blue}{"-----------------------------degree----------------------------", 3}

  \textcolor{blue}{"irreducible element of multidegree", [2, 1], "and order ",
        1, "found"}

  \textcolor{blue}{"irreducible element of multidegree", [0, 3], "and order ",
        3, "found"}

  \textcolor{blue}{"irreducible element of multidegree", [1, 2], "and order ",
        1, "found"}

  \textcolor{blue}{"-----------------------------degree------------------------------", 4}

  \textcolor{blue}{"irreducible element of multidegree", [2, 2], "and order ",
        0, "found"}

  \textcolor{blue}{"irreducible element of multidegree", [3, 1], "and order ",
        0, "found"}

  \textcolor{blue}{"irreducible element of multidegree", [0, 4], "and order ",
        0, "found"}

  \textcolor{blue}{"irreducible element of multidegree", [1, 3], "and order ",
        2, "found"}

  \textcolor{blue}{"-----------------------------degree----------------------------
        --", 5}

  \textcolor{blue}{"irreducible element of multidegree", [2, 3], "and order ",
        1, "found"}

  \textcolor{blue}{"-----------------------------degree----------------------------
        --", 6}

  \textcolor{blue}{"irreducible element of multidegree", [3, 3], "and order ",
        0, "found"}

  \textcolor{blue}{"-----------------------------degree----------------------------
        --", 7}

  \textcolor{blue}{"-----------------------------degree----------------------------
        --", 8}

  \textcolor{blue}{"-----------------------------degree----------------------------
        --", 9}

  \textcolor{blue}{"-----------------------------degree----------------------------
        --", 10}

  \textcolor{blue}{"-----------------------------degree----------------------------
        --", 11}

  \textcolor{blue}{"-----------------------------degree----------------------------
        --", 12}

  \textcolor{blue}{"-----------------------------degree----------------------------
        --", 13}

     \textcolor{blue}{"number  of semi-invariant in   minimal generating set ", 13}
\end{center}

 {\tt \textcolor{red}{> Ker[3];}}
 \textcolor{blue}{
 \begin{center}
\begin{multline*}
 \{ -6\,{y_{{0}}}^{2}y_{{3}}+6\,y_{{0}}y_{{1}}y_{{2}}-2\,{y_{{1}}}^{3},4\,x_{{1}}y_{{0}}y_{{2}}-6\,x_{{0}}y_{{0}}y_{{3}}+2\,x_{{0}}y_{{
1}}y_{{2}}-2\,x_{{1}}{y_{{1}}}^{2},-2\,{x_{{0}}}^{2}y_{{2}}-{x_{{1}}}^
{2}y_{{0}}+2\,x_{{0}}x_{{1}}y_{{1}} \} 
 \end{multline*}\end{center}
}
 
 \subsection{Compute $\mathcal{O}(V_1 \oplus V_1 \oplus V_1 \oplus V_2 \oplus V_2)^{SL_2}$ }  Use the command
 
 {\tt \textcolor{red}{ > dd:=[1,1,1,2,2]:Inv:=Min\_Gen\_Set\_Invariants\_S(dd):}}

\begin{center}

                \textcolor{blue}{  "calculating  Poincare series...."}

                       \textcolor{blue}{"done!, upper  bound", 6}

  \textcolor{blue}{"---------------------------------------------------------------
        ------------------"}

  \textcolor{blue}{"-----------------------------degree----------------------------
        --", 2}

                  \textcolor{blue}{6, "irreducible invariants  found"}

  \textcolor{blue}{"-----------------------------degree----------------------------
        --", 3}

                 \textcolor{blue}{12, "irreducible invariants  found"}

  \textcolor{blue}{"-----------------------------degree----------------------------
        --", 4}

                  \textcolor{blue}{" an irreducible invariant found"}

                  \textcolor{blue}{" an irreducible invariant found"}

                  \textcolor{blue}{" an irreducible invariant found"}

                  \textcolor{blue}{" an irreducible invariant found"}

                  \textcolor{blue}{" an irreducible invariant found"}

                  \textcolor{blue}{" an irreducible invariant found"}

  \textcolor{blue}{"-----------------------------degree----------------------------
        --", 5}

  \textcolor{blue}{"-----------------------------degree----------------------------
        --", 6}

                  \textcolor{blue}{"Total number of  invariants", 24}
\end{center}

{\tt \textcolor{red}{ > Ker[3]; }}
\textcolor{blue}{
\begin{multline*}
\{ 3\,{y_{{0}}}^{2}v_{{2}}-6\,y_{{0}}y_{{1}}v_{{1}}+3\,{y_{{1}}}^{2}v_{{0}},3\,{y_{{0}}}^{2}w_{{2}}+3\,{y_{{1}}}^{2}w_{{0}}-6\,y_{{0}}
y_{{1}}w_{{1}},6\,u_{{0}}u_{{1}}w_{{1}}-3\,{u_{{0}}}^{2}w_{{2}}-3\,{u_
{{1}}}^{2}w_{{0}},\\-6\,y_{{0}}u_{{0}}v_{{2}}-6\,y_{{1}}u_{{1}}v_{{0}}+6
\,y_{{1}}u_{{0}}v_{{1}}+6\,y_{{0}}u_{{1}}v_{{1}},-6\,y_{{0}}u_{{0}}w_{
{2}}-6\,y_{{1}}u_{{1}}w_{{0}}+6\,y_{{0}}u_{{1}}w_{{1}}+6\,y_{{1}}u_{{0
}}w_{{1}},\\6\,u_{{0}}u_{{1}}v_{{1}}-3\,{u_{{1}}}^{2}v_{{0}}-3\,{u_{{0}}
}^{2}v_{{2}},6\,u_{{1}}x_{{0}}v_{{1}}+6\,u_{{0}}x_{{1}}v_{{1}}-6\,u_{{0
}}x_{{0}}v_{{2}}-6\,u_{{1}}x_{{1}}v_{{0}},\\6\,y_{{1}}x_{{0}}w_{{1}}-6\,
y_{{1}}x_{{1}}w_{{0}}-6\,y_{{0}}x_{{0}}w_{{2}}+6\,y_{{0}}x_{{1}}w_{{1}
},-6\,u_{{1}}x_{{1}}w_{{0}}-6\,u_{{0}}x_{{0}}w_{{2}}+6\,u_{{0}}x_{{1}}
w_{{1}}+6\,u_{{1}}x_{{0}}w_{{1}},\\6\,y_{{0}}x_{{1}}v_{{1}}-6\,y_{{1}}x_
{{1}}v_{{0}}+6\,y_{{1}}x_{{0}}v_{{1}}-6\,y_{{0}}x_{{0}}v_{{2}},-3\,{x_
{{1}}}^{2}v_{{0}}+6\,x_{{0}}x_{{1}}v_{{1}}-3\,{x_{{0}}}^{2}v_{{2}},-3
\,{x_{{0}}}^{2}w_{{2}}+6\,x_{{0}}x_{{1}}w_{{1}}-3\,{x_{{1}}}^{2}w_{{0}
} \}
 \end{multline*}
 }
\end{document}